\documentclass[12pt,openany,twoside]{memoir}

\usepackage{newtxtext}
\usepackage[T1]{fontenc}
\usepackage{mathpazo}

\usepackage{amsmath,amsthm,amssymb}
\usepackage{hyperref}

\usepackage{aliascnt}
\usepackage{xcolor}

\definecolor{dark}{rgb}{0,0,0.6}

\definecolor{light-gray}{gray}{0.7}

\hypersetup{
	unicode=true,
	colorlinks=true,
	citecolor=dark,
	linkcolor=dark,
	anchorcolor=dark
}

\DeclareTextFontCommand{\citefont}{\fontfamily{qtm}\mdseries}

\usepackage{cite}

\settrims{0pt}{0pt}
\settypeblocksize{22.5cm}{17cm}{*}
\setmarginnotes{10pt}{10pt}{10pt}
\setbinding{.5cm}
\setulmargins{*}{*}{*}
\setlrmargins{*}{*}{*}
\checkandfixthelayout

\makechapterstyle{myell}{%
  \chapterstyle{default}
  \renewcommand*{\chapnumfont}
  				{\normalfont\color{light-gray}\fontsize{58}{58} \fontfamily{qtm}\mdseries}
  
  \settowidth{\chapindent}{\chapnumfont 55}
  \setlength{\chapindent}{3em}

  }

\chapterstyle{myell}

\usepackage{fancyhdr}

\fancypagestyle{plain}{%
	\fancyhf{} 
	\fancyfoot[C]{\normalfont \thepage} 

}

\checkandfixthelayout

\pretitle{\begin{center}\Huge\bfseries}
\posttitle{\par\end{center}}
\newcommand{\submitted}[1]{\gdef\suB{#1}}
\newcommand{\revised}[1]{\gdef\puB{#1}}
\newcommand{\suB}{}
\newcommand{\puB}{}

\newtheorem{theorem}{Theorem}

\newaliascnt{lemma}{theorem}
\newaliascnt{proposition}{theorem}
\newaliascnt{corollary}{theorem}
\newaliascnt{fact}{theorem}
\newaliascnt{example}{theorem}

\newtheorem{lemma}[lemma]{Lemma}
\newtheorem{proposition}[proposition]{Proposition}
\newtheorem{corollary}[corollary]{Corollary}

\aliascntresetthe{lemma}
\aliascntresetthe{proposition}
\aliascntresetthe{corollary}
\aliascntresetthe{fact}
\aliascntresetthe{example}

\theoremstyle{definition}
\newtheorem{definition}[theorem]{Definition}
\newtheorem*{remark}{Remark}
\newtheorem*{question}{Question}

\newcommand{\axiom}[1]{\mathsf{#1}}
\newcommand{\ZF}{\axiom{ZF}}
\newcommand{\ZFC}{\axiom{ZFC}}
\newcommand{\DC}{\axiom{DC}}
\newcommand{\AC}{\axiom{AC}}

\newcommand{\BPI}{\axiom{BPI}}
\newcommand{\BP}{\axiom{BP}}

\newcommand{\RR}{\mathbb R}
\newcommand{\CC}{\mathbb C}
\newcommand{\NN}{\mathbb N}
\newcommand{\QQ}{\mathbb Q}
\newcommand{\ZZ}{\mathbb Z}
\newcommand{\cF}{\mathcal F}

\newcommand{\power}{\mathcal P}

\newcommand{\tup}[1]{\langle#1\rangle}

\DeclareMathOperator{\rng}{rng}
\DeclareMathOperator{\spnn}{span}
\DeclareMathOperator{\dom}{dom}

\renewcommand{\d}{\,\mathrm{d}}

\begin{document}
\frontmatter
\pagenumbering{gobble}
\title{Zornian Functional Analysis\\\par\large{or: How I Learned to Stop Worrying and Love the Axiom of Choice}}
\author{Asaf Karagila}
\submitted{February 1, 2016}
\revised{October 29, 2020}
\date{}
\maketitle
\newpage
\pagenumbering{roman}
\pagestyle{fancy}
\setcounter{page}{1}
\tableofcontents
\chapter{Preface}
This work is the aftereffect of a course about the fundamentals of functional analysis given by Professor Matania Ben-Artzi during the fall semester of 2015--2016. The title references Garnir's paper \cite{Garnir:1973}, as well Stanley Kubrick's masterpiece ``\textit{Dr.\ Strangelove or: How I Learned to Stop Worrying and Love the Bomb}''. The idea is that one should embrace the axiom of choice, rather than reject it. This work might seem a bit strange in this context, as it investigates what might happen without the presence of the axiom of choice, but one should view it as a motivation for accepting the axiom of choice. One could say that this is a proof-by-intimidation (if you will not assume the axiom of choice such and such will happen!) and the author will plead no contest to such accusations.

Much of this knowledge was obtained by the author by digging up and gathering papers over the years preceding the course, and this was the perfect time to sit and review many of these proofs with a matured eye and slightly better understanding of functional analysis. In all likelihood the majority of these proofs can be found in the ``\textit{Handbook of Analysis and its Foundations}'' by Eric Schechter \cite{Schechter:1997}, which can be an excellent source for many results about the foundations of analysis and particularly in relation to the axiom of choice. Another accessible source is ``\textit{Axiom of Choice}'' by Horst Herrlich \cite{Herrlich:AC} which is a book full of references and information accessible to the working mathematician about what can happen without the axiom of choice.

The intended reader is not your average set theorist (or even your average set theory graduate students). The ideal reader is rather someone who understands the basics the functional analysis, and we hope to have achieved a reasonable level of clarity for that crowd. Knowledge in set theory, while always positive, is not really needed if you are willing to believe the author.

Lastly, it is important to point out that this work contains no original research, although many of the proofs presented here have been digested and presented in the image of the author (and any mistakes are likely his fault).

\subsubsection*{Acknowledgments.} This work would not have been written (at least not in its present form) without the kind permission of Prof.\ Ben-Artzi to submit this work as a final project; Theo B\"uhler whose encyclopedic knowledge regarding historical results and various references is uncanny, and helped in clearing out some issues; Matt Foreman who suggested to include the part about Banach limits; and finally Ariel Rapaport for his availability during the writing of this work for small consultations.

\mainmatter
\chapter{Introduction}
This text is meant for analysis students who want to learn more about the effects of the axiom of choice on functional analysis, and the things that may go wrong in its absence. As this is a text aimed for analysis students, we will not focus on the set theoretic proofs or dwell on any particular set theoretic assumptions that are needed in order to prove certain results.

For those interested in reading more about the set theoretic side of these constructions we will provide ample references and citations, we recommend Enderton's book ``\textit{Elements of Set Theory}'' \cite{Enderton:Elements} as introductory to set theory, then continuing with Jech \cite{Jech:ST}, Kunen \cite{Kunen:2011Book} or Halbeisen \cite{Halbeisen:CombST} (the last is particularly suitable for choice related independence proofs).

From a slightly more accurate foundational point of view, this work is written within the context of classical first-order logic and the set theoretic axioms of Zermelo--Fraenkel (abbreviated as $\ZF$ and $\ZFC$ when the axiom of choice is adjoined to our theory). We will not get into the axioms here, but naively they tell us that sets behave in the way we expect sets to behave.

\paragraph{What do we assume to be known by the reader?} This text is not aimed at people who have just started their first analysis course, nor their second. The assumption is that the reader is fluent in the basics of topology, linear algebra, measure theory, and has at least some familiarity with the functional analysis covered here. From the set theoretic aspects, we \textbf{do not assume} the reader is fluent in the language of ordinals, cardinals and so on. However we do expect rudimentary knowledge of set theory which includes (but not limited to) the basics of set theory, partial orders, and so on. Hopefully, after this text the reader will want to abandon their plans in order to pursue a career in set theory.

\newpage
\section{What is the axiom of choice?}
\begin{definition}[\textit{The Axiom of Choice}]\label{Def:AC}
If $\{A_i\mid i\in I\}$ is a set of non-empty sets, then there exists a function $f$ with domain $I$ such that $f(i)\in A_i$ for all $i\in I$.
\end{definition}

In simple words, the axiom of choice asserts that if we have a collection of non-empty sets, then we can ``choose'' an element from each one. Mathematically speaking, this means that there is a function which takes in $A_i$ (or its index), and gives us an element which belongs to $A_i$. Such a function is called a \textit{choice function}.

\subsection{Some historical background}
When Cantor started developing set theory, one of the theorems he proved was that given two sets $A$ and $B$, there is an injection from $A$ to $B$ or vice versa. His proof relied on the fact that every set can be well-ordered, which was taken as a naive truth. In 1904 Zermelo proposed an axiomatic treatment of set theory, and included the axiom of choice in order to prove that every set can be well-ordered.

The other axioms of set theory (those originally given by Zermelo, and the later modifications by Fraenkel and Skolem into what is now known as $\ZF$) are axioms which assert the existence of specific sets. For example, if $A$ is a set, then its power set exists by the axiom of power sets, and we know exactly what set is $\mathcal P(A)$; if $\varphi(x)$ is some property and $A$ is a set, then by the subset axiom, $\{a\in A\mid \varphi(a)\}$ is a set as well, so if you know what is $\varphi$, then you know the resulting set as well. The axiom of choice, however, asserts only the existence of a choice function and it tells us nothing  else.

Due to its ``mysterious'' and non-constructive nature, many mathematicians initially rejected the axiom and developed what eventually culminated in the constructive school of mathematics and intuitionistic logic. When one takes into consideration many of the counter-intuitive consequences (like the existence of non-measurable sets) which follow the axiom of choice, it is not surprising at all that dissenting voices were heard.\footnote{Gregory Moore wrote in details on the history of the axiom of choice in \cite{Moore:AC}.}

However, the axiom of choice allows us to bring structure and order into the world of infinite sets. Its usefulness triumphed, and after G\"odel proved in 1938\footnote{See \cite{Godel:AC}, or any book about modern set theory, e.g.\ \cite{Kunen:2011Book}.} that the axiom of choice does not add new contradictions to $\ZF$, mathematicians accepted it into the standard assumptions of mathematics. In 1963, when Cohen proved \cite{Cohen:1963, Cohen:1964} that the axiom of choice is indeed not provable from $\ZF$, the axiom of choice became a topic of investigation by set theorists for the past 50 years (see \cite{HowardRubin:Consequences, Jech:AC} for examples).

The set theoretic research about the axiom of choice involves various fragments of the axiom of choice, which arise naturally in various proofs throughout mathematics. The research is mainly focused in the interactions between the various fragments of the axiom, and the interaction with ``common'' mathematical results (e.g.\ Baire Category Theorem). The latter will be the focus of this work, in relation to analysis.

\subsubsection{The relation to functional analysis}
The axiom of choice was always linked to functional analysis. From Banach onward the axiom of choice was used for proving theorems like Hahn--Banach theorem, Baire's theorem, and even before that sequential continuity was proved to be equal to $\varepsilon$-$\delta$ continuity. It was after Cohen's discovery of forcing and the availability of tools to produce models of set theory in which the axiom of choice failed that the necessity of the axiom was discovered.

In the late 1950s and early 1960s model theoretic results like ultraproducts were used to deduce from model theoretic choice principles (e.g.\ the compactness theorem) statements like the Hahn--Banach theorem, notable works are \L{}os and Ryll-Nardewski \cite{Los-RyllN}, and Luxemburg \cite{Lux:1962}. At that point it was already known that the various model theoretic techniques rely heavily on the axiom of choice.

The 1970s were rife with discovery. In 1970 Solovay published his seminal result that $\ZF$\footnote{Along with a weak form of the axiom of choice, sufficient for all the basic measure theoretic definitions and theorems to be proved as usual.} is consistent with the assumption that every set of reals was Lebesgue measurable, and more \cite{Solovay:1970}. It quickly turned out that in Solovay's model functional analysis was completely changed. Garnir in his paper \cite{Garnir:1973} suggested that functional analysis will be split into Zornian, Solovayan and Constructive variants depending on which side of the axiom of choice you prefer. Garnir's paper, as well other papers by Wright \cite{Wright:1973} and V\"ath \cite{Vath:1998}, are entirely devoid of set theory and they focus solely on the functional analysis that can be done under the axioms which hold in Solovay's model.

Pincus and others \cite{Pincus-Solovay:1977, Pincus:1974} made other related refinements to explore the logical strength of various functional analytic theorems. Bell and Fremlin \cite{Bell-Fremlin:1972} proved that Hahn--Banach in conjunction with a mild strengthening of Krein--Milman theorem can be used to prove the axiom of choice, for example.

Alas, Garnir's dream mathematics required functional analysts to give up some of the most useful tools in their toolbox, e.g.\ Hahn--Banach, and the majority of functional analysis remained within the domain of $\ZFC$. But even nowadays there is active research into the necessity of the axiom of choice in functional analysis, and many questions remain open. It should be added, perhaps, that while many mathematicians regard the axiom of choice as an obvious truth, there is still importance to this research. It tells us exactly what sort of things in our work cannot be explicitly constructed. If functional analysis eventually trickles down to physics and engineering this understanding is something of philosophical and practical significance.

\section{Some basic results related to the axiom of choice}
In this section we will present, without proof, a few equivalents and consequences of the axiom of choice. For a far more complete list we encourage the reader to look at the books \cite{HowardRubin:Consequences,RubinRubin:Equivalents} which are dedicated to statements equivalents to the axiom of choice from different fields of mathematics.\newpage
\subsection{Common equivalents}
The two most famous equivalents of the axiom of choice are Zorn's lemma and the well-ordering principle. But there are a few more useful equivalents.
\begin{definition}\label{Def:WO}
We say that $\tup{A,<}$ is a \textit{well-ordered set} if $<$ is a linear ordering of $A$ and whenever $B\subseteq A$ is a non-empty set, then $B$ has a $<$-minimal element.
\end{definition}
The idea behind well-orders is that they extend the natural numbers in a way that allows us to perform induction and recursion over the linear order. Such induction is called \textit{transfinite induction} (or transfinite recursion).

\begin{theorem}
The following are equivalent:
\begin{enumerate}
\item The axiom of choice.
\item Every set can be well-ordered.
\item Zorn's lemma: If $\tup{P,<}$ is a partially ordered set in which every chain has an upper bound, then there is a maximal element in $P$.
\item Hausdorff's Maximality Principle: If $\tup{P,<}$ is a partially ordered set, then there is a maximal chain in $P$.
\item If $f\colon A\to B$ is a surjective function, then there is an injective $g\colon B\to A$ such that $f(g(b))=b$ for all $b\in B$.
\item Given two sets $A$ and $B$, there is either an injection $f\colon A\to B$ or an injection $f\colon B\to A$.
\item Every vector space has a basis.
\end{enumerate}
\end{theorem}

We can see why the axiom of choice is a popular assumption. It helps brings order and structure to the world of infinite sets. And as mathematics evolved beyond the world of finite objects, into the realm of infinite sets, we needed more assumptions in order to make sense of that world. However, this order stretches well beyond the interest of the common mathematician. So what sort of weaker versions of the axiom of choice are sufficient for most wants of the common mathematician?

\subsection{Weak choice axioms}

The most common restriction of the axiom of choice is perhaps the axiom of countable choice. This axiom asserts that given a countable family of non-empty sets, then we can find a choice function for it. However in practice this is not the most useful choice principle. One of the most implicit uses of the axiom of choice is when defining a sequence by recursion and making arbitrary choices along the way. It is a common mistake that the axiom of countable choice is sufficient for these sort of constructions, but in fact it is not.

\begin{definition}The principle of \textit{Dependent Choice} asserts that given a non-empty set $X$ and a binary relation $R\subseteq X\times X$, such that for every $x\in X$ there is some $y\in X$ for which $x\mathrel{R}y$, then for every $x_0\in X$ there is a function $f\colon\NN\to X$ such that $f(0)=x_0$ and for every $n\in\NN$, $f(n)\mathrel{R} f(n+1)$. We will often abbreviate this statement as $\DC$.
\end{definition}

If one stops to ponder this definition, it is not hard to see that it says exactly that recursive definitions work. Given a set $X$ if we want to define a sequence by recursion this amounts to picking some $x_0$, and proving that given $x_n$ there is a non-empty selection for viable candidates to serve as $x_{n+1}$. The strength of Dependent Choice is nontrivial, but it is still quite far from the entire axiom of choice.

\begin{theorem}
Each of these statements imply the next, and no implication can be reversed:\footnote{The proofs can be found in \cite{Jech:AC}.}
\begin{enumerate}
\item The axiom of choice.
\item The principle of dependent choice.
\item The axiom of countable choice.
\item The countable union of countable sets is countable.
\end{enumerate}
\end{theorem}

These implications are clearly very useful in our everyday mathematics, and along with the reinterpretation of $\DC$ as allowing us generalized recursive definitions, it seems that $\DC$ is a very handy tool for measure theory and basic functional analysis. However, $\DC$ is not the most useful to the working mathematician. This title belongs to the principle often called the Boolean Prime Ideal theorem, or $\BPI$. Presenting $\BPI$ properly requires introducing some definitions that will not be of much use in this work, so instead we present an equivalent statement.

\begin{definition}
The \textit{Boolean Prime Ideal theorem} asserts that for every set $X$ and a collection $\mathcal D$ of subsets of $X$ with the finite intersection property\footnote{Any finite subcollection has a non-empty intersection} there is a two-valued measure $\mu$ on $\power(X)$ such that $\mu(A)>0$ for all $A\in\mathcal D$.\footnote{In here the two-valued measure can be replaced by an ultrafilter; or by requiring that every filter can be extended to an ultrafilter; we may also replace $\power(X)$ with an arbitrary Boolean algebra and simply require the existence of an ultrafilter.} We will abbreviate this statement as $\BPI$. However abstract Boolean algebras and filters are exceeding the scope of this work, so we will focus on two-valued measures.
\end{definition}

To establish the claim that $\BPI$ is one of the most of useful choice principles in mathematics, we will list a few equivalents and a few consequences. This is by no means a complete list.

\begin{theorem}
The following are equivalent:\footnote{Proofs can be found in \cite{Jech:AC,Herrlich:AC}.\label{foot7}}
\begin{enumerate}
\item $\BPI$.
\item If $\{X_i\mid i\in I\}$ is a family of compact Hausdorff spaces,\footnote{We do not follow Bourbaki's convention that compact implies Hausdorff.} then $\prod_{i\in I}X_i$ is a compact Hausdorff space.
\item For every $I$, $\{0,1\}^I$ is compact, where $\{0,1\}$ is the discrete space with two elements.
\item The compactness theorem for first-order logic.
\item The completeness theorem for first-order logic.
\item If $R$ is a commutative ring with a unit, then every ideal is contained in a prime ideal.
\item Stone--\v{C}ech compactification theorem.
\item Stone representation theorem.
\item If $X$ is a Banach space, then the closed unit ball of $X^*$ is weak-$*$ compact.
\end{enumerate}
\end{theorem}

\begin{theorem}
The following are consequences of $\BPI$. None of these statements is provable in $\ZF$ alone, and none of them imply $\BPI$.\footref{foot7}
\begin{enumerate}
\item Hahn--Banach theorem.
\item Every set can be linearly ordered.
\item There exists a non-measurable set.
\item The Banach--Tarski paradox.
\item There is a finitely additive probability measure on $\power(\NN)$ which vanishes on singletons.
\item If a vector space has a basis, then every two bases have the same cardinality.
\item Every two algebraic closures of a field $F$ are isomorphic.
\end{enumerate}
\end{theorem}

In the coming chapters we will mainly work in one of the following theories: $\ZF$, $\ZF+\DC$, $\ZF+\BPI$ or $\ZF+\DC+\BPI$. Any assumptions exceeding $\ZF$ will be explicitly stated in the theorems. But before we do that, we also want to exhibit how badly things can fail in the absence of choice.

\subsection{Failures of interest}

The axiom of choice is often criticized as being non-constructive, as it asserts the existence of certain objects, but tells us nothing as to how they ``look like''. The problem with the failure of the axiom of choice is that it is just as non-constructive. If the axiom of choice fails, we have no way knowing a priori what family of sets does not admit a choice function, or what vector space has no basis, and so on. It is true, that from a counterexample to one formulation or equivalent, we can usually manufacture relatively specific counterexamples to other formulations and equivalents to the axiom of choice. It does not make the original counterexample which we can only prove to exist any more constructive.

If we just assume the failure of the axiom of choice, it also does not tell us whether or not slightly weaker assertions still hold, or if it even fails at any level of interest to the non-set theorist. For example it might be that the axiom of choice fails in several acute ways, but these failures happen somewhere far away from the real numbers or any vector space, ring or module which interest the majority of mathematicians. It is consistent, though, that the axiom of choice fails in ways so severe that it has implications to the common mathematics. We will list a few of these in order to motivate some basic assumptions later on.
\begin{theorem}
The following statements are consistent with the failure of the axiom of choice:
\begin{enumerate}
\item $\RR$ is the countable union of countable sets.\footnote{Note that $\RR$ is provably uncountable even without the axiom of choice.}
\item There exists a metric space $X$ and a function from $X$ into $\RR$ which is sequentially continuous, but not $\varepsilon$-$\delta$ continuous.
\item $\NN$ is not a Lindel\"of space, equivalently $\RR$ is not a Lindel\"of space.
\item There is a metric space which is not compact, but it is sequentially compact.
\item There is a dense subset of $\RR$ in the Euclidean topology which does not have a countably infinite subset.
\item $\RR$ is the union of two sets, each of a strictly smaller cardinality.
\item The only linear functional on $\ell^\infty/c_0$ is the $0$ functional.
\item $(\ell^\infty)^*=\ell^1$.
\item Every linear functional on a Fr\'echet space is continuous.
\item Every set of real numbers is Lebesgue measurable.
\item Every set of real numbers is Baire measurable.
\item $\QQ$ has two non-isomorphic algebraic closures.
\end{enumerate}
\end{theorem}

It is worth to mention that $\BPI$ is consistent with the statement ``There is an infinite set of real numbers which does not have a countably infinite subset''. The latter is a contradiction to $\DC$ and therefore we cannot prove $\DC$ in $\ZF+\BPI$. In the other direction it is consistent with $\ZF+\DC$ that every set is Lebesgue measurable, in which case $\BPI$ fails and therefore we cannot prove $\BPI$ from $\ZF+\DC$. So these two principles are indeed independent. Moreover, while $\ZF+\DC+\BPI$ seems to be sufficient to prove a lot of results of interest, this theory still cannot prove the axiom of choice in full.

We conclude this section by pointing out that in a model where $\RR$ is the countable union of countable sets, we cannot begin to talk about measure theory, since the only nontrivial $\sigma$-additive measures on $\RR$ (or any set of equal cardinality) are atomic. Why? If every singleton has measure $0$, then every countable set has measure $0$, and therefore every countable union of countable sets has measure $0$. Unfortunately, this includes the entire space. If one was not convinced before, this makes a grand case for accepting $\DC$ as a valid mathematical axiom.

\chapter{Baire Category Theorem}

\section{Baire's Theorem}
\begin{definition}
Let $X$ be a topological space and $A\subseteq X$.
\begin{enumerate}
\item We say that $A$ is \textit{dense} if $\overline{A}=X$, or equivalently if for every non-empty open set $U$ it holds that $A\cap U\neq\varnothing$.
\item We say that $A$ is \textit{nowhere dense} if $X\setminus\overline A$ is dense.
\item We say that $A$ is \textit{meager} or \textit{first category} if it is the countable union of closed nowhere dense sets.
\item We say that $A$ is \textit{co-meager} if $X\setminus A$ is meager. Note that if $X$ is not meager, then being co-meager is stronger than not being meager.
\end{enumerate}
\end{definition}

\begin{definition}
Let $X$ be a topological space. We say that $X$ is a \textit{Baire space} if the countable intersection of dense open sets is dense in $X$.
\end{definition}

\begin{theorem}
Let $X$ be a Baire space, then $X$ is not meager.
\end{theorem}
\begin{proof}
Suppose that $E_n$ is a closed nowhere dense set for $n\in\NN$, then $D_n=X\setminus E_n$ is a dense open set. Therefore $D=\bigcap_{n\in\NN}D_n$ is dense, and in particular non-empty. Therefore $\bigcup_{n\in\NN}E_n\neq X$.
\end{proof}

The eponymous theorem of Baire is the following theorem.

\begin{theorem}[$\AC$]
If $X$ is a complete metric space, then $X$ is a Baire space.
\end{theorem}
\begin{proof}
For $n\in\NN$ let $D_n$ be a dense open subset of $X$. We may assume that $D_{n+1}\subseteq D_n$ for all $n$, since the finite intersection of dense open sets is a dense open set again, so we can replace $D_n$ by $\bigcap_{k=1}^n D_n$ if needed.

Let $\varnothing\neq U\subseteq X$ be an open set. There is some $r_0>0$ and $x_0\in U$ such that $B(x_0,r_0)\subseteq U\cap D_0$. Recursively define $x_n,r_n$ such that $0<r_{n+1}<2^{-n}r_n$ and $B(x_{n+1},r_{n+1})\subseteq B(x_n,r_n)\cap D_n$. It is clear that $\{x_n\}_{n=1}^\infty$ is a Cauchy sequence and by completeness it has a limit $x$ which has to be an element of each $B(x_n,r_n)$ and therefore $x\in\bigcap_{n\in\NN}D_n\cap U$, so the intersection is indeed dense.
\end{proof}

The reader should perhaps feel a bit unease seeing the above proof assumes $\AC$ rather than just $\DC$. And indeed we only used $\DC$ here. But more on that later.

\section{Disasters without choice}
For example, in a $T_1$ space, where every singleton is closed, every countable set is meager. We mentioned that it is consistent that $\RR$ is the countable union of countable sets, but does that mean that $\RR$ itself is meager in that situation? What other disasters can happen if the axiom of choice is not assumed?

First of all, when we examine the proof the only place the axiom of choice is used is when defining the sequences of $x_n$'s and $r_n$'s, because there is really no reason to prefer one point and one radius over the other. We can easily dispose of the choice used to pick the $r_n$'s by limiting ourselves to radii of the form $\frac1n$ for some $n\in\NN$, in which case we can always take the least possible $n$ to satisfy the inequality.

But picking $x_n$ canonically? There is no reason that we should be able to do that. And indeed we might not be able to do that.

Recall now that in our failures of interest we remarked it is consistent that there is $X\subseteq\RR$ which is dense, and thus infinite, but does not have a countably infinite subset.

Define the metric on $X^\NN$ (all the functions from $\NN$ to $X$) as usual on product spaces,\[ d(f,g)=2^{-\min\{n\mid f(n)\neq g(n)\}}.\]
\begin{theorem}
$X^\NN$ is a complete metric space with $d$ as defined above.\qed
\end{theorem}

\begin{theorem}\label{thm:failure of baire}
The space $X^\NN$ is not a Baire space.
\end{theorem}
\begin{proof}
For every $n\in\NN$ define the set $U_n$ to be
\[U_n=\bigcup_{m=n+1}^\infty\{f\in X^\NN\mid f(n)<f(m)\}\]
(recall that $X$ is a set of real numbers, so $f(n)<f(m)$ is a meaningful statement). We first claim that $U_n$ is dense open for every $n\in\NN$.

To see that $U_n$ is open, suppose that $f\in U_n$, then there is some $m$ such that $f(n)<f(m)$. Now if $d(f,g)<2^{-m-1}$ it follows that $f(i)=g(i)$ for all $i<m+1$ and in particular $g(n)=f(n)$ and $g(m)=f(m)$ so $g\in U_n$. The density follows from the fact that if $f\in X^\NN$ and $m>n$, fix some $x\in X$ such $f(n)<x$, then define $g\colon\NN\to X$, \[g(k)=\begin{cases} f(k) & k\neq m\\ x & k=m\end{cases}\]
It follows that $g\in U_n$ and $d(f,g)=2^{-m}$, and therefore we can find points in $U_n$ which are arbitrarily close to $f$.

Now we claim that $\bigcap U_n$ is not dense. In fact, we claim that it is empty. If $f\in U_n$ for all $n\in\NN$, then for every $n$ there is $m>n$ such that $f(n)<f(m)$; without using the axiom of choice we can find an infinite subset of $\NN$ on which $f$ is injective. This is a contradiction as $X$ does not have a countably infinite subset, and the range of $f$ is countably infinite.
\end{proof}

The above proof will be helpful later, but we will remark that in the original construction of a set theoretic universe in which there is a subset of $\RR$ which is dense but has no countably infinite set, such subset will be a complete metric space but not a Baire space. The proof of completeness is trivial, since a Cauchy sequence must be eventually constant (otherwise it is either not Cauchy, or defines a countably infinite subset), so as a metric space with the metric induced from $\RR$ it is complete. The proof that this space is not a Baire space requires looking into the bowels of the set theoretic constructions and arguments, and will not be brought here.

\section{Do we really need the axiom of choice?}

In the previous section we saw a counterexample to the classic theorem that every complete metric space is a Baire space. We used the failure of the axiom of choice in which there was an infinite set of real numbers without a countably infinite subset. consider the following observation:
\begin{theorem}[$\DC$]
Let $X$ be a topological space. The countable union of meager sets is meager.
\end{theorem}
\begin{proof}
Let $M_n$ be meager sets for $n\in\NN$ and $E_{n,k}$ for $k\in\NN$ closed nowhere dense sets such that $M_n=\bigcup_{k\in\NN}E_{n,k}$. Then $\bigcup_{n\in\NN}M_n=\bigcup_{(n,m)\in\NN\times\NN}E_{n,k}$ is the countable union of closed nowhere dense sets and therefore meager.
\end{proof}

The above proof is a bit misleading, as it seems not to be using the axiom of choice at all (or even $\DC$). The axiom is used in choosing a sequence of closed sets for each $M_n$ (along with its enumeration).We also know that the countable union of countable sets need not be countable in $\ZF$, e.g.\ $\RR$ can be the countable union of countable sets. Certainly countable sets are meager in $\RR$, does that mean that $\RR$ is meager in such situation? Certainly this would be a far more spectacular counterexample to Baire's theorem. Fortunately, this is not possible, as shown by the following theorem.

\begin{theorem}
If $X$ is a complete, separable metric space, then $X$ is a Baire space. In particular $\RR$ is a Baire space.
\end{theorem}
\begin{proof}
Fix $\{e_n\mid n\in\NN\}$ a countable dense subset of $X$. Now let $U$ be a non-empty open set and $D_n$ a sequence of dense open sets. Repeat the proof of Baire's theorem, only now let $x_{n+1}$ be $e_{k(n)}$ where $k(n)$ is the least $k\in\NN$ such that $d_k\in B(x_n,r_n)\cap D_n$, and let $r_{n+1}$ be the least $m\in\NN$ such that $\frac1m<2^{-n}r_n$ and $B(x_{n+1},r_{n+1})\subseteq B(x_n,r_n)\cap D_n$.

Again $\{x_n\}_{n=1}^\infty$ is a Cauchy sequence so it must converge to some $x$ in $\bigcap_{n\in\NN}D_n\cap U$.
\end{proof}

This theorem, as remarked, guarantees us Baire's theorem for $\RR$ and other important spaces like the finite dimensional $\RR^n$ and $\CC^n$ as well $L^p$ spaces\footnote{At least for $\sigma$-finite measure spaces.} for $1\leq p<\infty$. All those will remain Baire spaces in the usual applications even if the axiom of choice is not assumed at all.

\section{On consequences of the Baire Category Theorem}

\begin{theorem}[Blair \cite{Blair:1977}]
Baire Category Theorem is equivalent to $\DC$.
\end{theorem}
\begin{proof}
If we assume $\DC$ then the usual proof goes as usual. The interesting direction is to prove that $\DC$ is true, assuming Baire's theorem. Suppose that $R$ is a binary relation on $X$ such that for every $x\in X$ there is some $y\in X$ such that $x\mathrel{R}y$.

Define the metric on $X^\NN$ as in \autoref{thm:failure of baire}, this is a complete metric space and therefore a Baire space. For every $n\in\NN$ define as in the proof of \autoref{thm:failure of baire} the dense open set $U_n$,
\[U_n=\bigcup_{m=n+1}^\infty\{f\in X^\NN\mid f(n)\mathrel{R} f(m)\}.\]
As before $U_n$ is indeed dense open, so $\bigcap_{n\in\NN}U_n$ is dense, so it is not empty. Let $f\in U_n$ for all $n\in\NN$, then we have that $f\in U_n$, as in the proof of \autoref{thm:failure of baire} we get an infinite subset of $\NN$ on which $f$ is strictly ``increasing'', without loss of generality we can assume this to be $\NN$ itself, and so we get $f(n)\mathrel{R}f(n+1)$ for all $n\in\NN$ as wanted.
\end{proof}

In functional analysis Baire's theorem is used to prove the Open Mapping Theorem, Uniform Boundedness Theorem and the Closed Graph Theorem. These proofs do not use any more choice than Baire's theorem provides already. We can also easily prove these three consequences are equivalent in $\ZF$, so assuming one provides us with the others.

In recent events, Adrian Fellhauer showed that the Open Mapping Theorem can be proved using only the axiom of choice for countable families of sets. This assumption is strictly weaker than $\DC$, and so we get the following theorem.

\begin{theorem}[Fellhauer \cite{Fellhauer:2015}]
The Open Mapping Theorem does not imply Baire's theorem.\qed
\end{theorem}

In $\ZFC$ the two statements are provable, and therefore equivalent. But the source of their equivalence turns out to be the axiom of choice. If one adopts a slightly more constructive approach to philosophy of mathematics, then it means exactly that Baire's theorem is a stronger theorem than the Open Mapping Theorem.

\chapter{Dream Mathematics: Automatic Continuity}
We say that a topological vector space $X$ is a \textit{dream space} if every linear functional on $X$ is continuous. Every finite dimensional $\RR^n$ or $\CC^n$ is a dream space. Brunner investigated\footnote{And coined the term ``dream space''.} characterizations of dream spaces which have Hamel bases \cite{Brunner:1987}. In this chapter we will investigate the situation where every Fr\'echet space is a dream space.
\section{The Baire Property}
\begin{definition}Let $X$ be a topological space, we say that $A\subseteq X$ has the \textit{Baire property} if there exists an open set $U$ such that symmetric difference $A\mathop{\triangle} U$ is meager.
\end{definition}
Assuming the axiom of choice, a countable union of sets with the Baire property also has the Baire property. And in well-behaved spaces, like separable metric spaces, the sets which have the Baire property form a $\sigma$-algebra which contains all the open sets and so all the Borel sets. We will see later that the axiom of choice ($\BPI$ suffices) implies the existence of a set without the Baire property.

\begin{definition}
Let $X$ and $Y$ be topological spaces a function $f\colon X\to Y$ is \textit{Baire measurable} if $f^{-1}(U)$ has the Baire property in $X$ whenever $U$ is open in $Y$.
\end{definition}

It is easy to see that every continuous function is Baire measurable, but in the presence of choice, where in $\RR^n$ every Borel set has the property of Baire, we also have that every Borel measurable function is Baire measurable.

\subsection{Cauchy's functional equation}
We begin with this freshman calculus exercise.
\begin{theorem}
Suppose that $f\colon\RR\to\RR$ such that $f(x+y)=f(x)+f(y)$. If $f$ is continuous then $f(x)=ax$ for some $a\in\RR$.
\end{theorem}
\begin{proof}
Let $a=f(1)$, then by induction $f(n)=an$ for every $n\in\NN$ and therefore also for $n\in\ZZ$. Next, for every $\frac pq\in\QQ$ with $q>0$, we have that $\sum_{i=1}^q f(\frac pq)=f(p)=a p$ and therefore $f(\frac pq)=a\frac{p}q$. So for every rational number it holds that $f(x)=ax$. Finally, if $x\in\RR$, pick some $\{q_n\}_{n=1}^\infty\subseteq\QQ$ such that $q_n\to x$, then $f(q_n)\to f(x)$ by continuity and therefore $f(x)=ax$ as well.
\end{proof}

What if we weaken the requirement from being continuous to being Baire measurable? The answer turns out to be the same. We will not prove this theorem now, instead we will discuss a far more general theory with far greater consequences.

\section{Every set of reals has the Baire property}

\begin{theorem}[Solovay \cite{Solovay:1970}, Shelah \cite{Shelah:1984}\footnote{The result was first proved by Solovay in 1970 under assumptions that exceed $\ZFC$ that it is consistent that every set of reals is Lebesgue measurable and has the Baire property. In 1984 Shelah proved that for the Lebesgue measurability the added assumptions are necessary, but if we are only interested in the Baire property then these additional assumptions can be waived.}]
It is consistent with $\ZF+\DC$ that every set of reals has the Baire property.\qed
\end{theorem}
We shall abbreviate the statement that every set of reals has the Baire property as $\BP$. The following is a corollary of the fact that every uncountable complete and separable metric space has the same Borel and Baire structure as the reals.
\begin{corollary}
If every set of reals has the Baire property, then for every uncountable complete and separable metric space, every set has the Baire property as well.\qed
\end{corollary}

\subsection{Automatic continuity of linear operators}\label{Subsection:Auto}
\begin{lemma}[$\AC$]
Suppose that $X$ is a separable complete metric space and $Y$ is a separable metric space. If $f\colon X\to Y$ is Baire measurable, then there is a meager set $N$ such that $f\restriction(X\setminus N)$ is continuous.
\end{lemma}
\begin{proof}
Let $\{y_n\mid n\in\NN\}$ be a countable dense subset of $Y$ and consider $S_n=B_Y(y_n,\varepsilon/2)$ for some fixed $\varepsilon>0$. Note that $\{S_n\mid n\in\NN\}$ is a cover of $Y$. Next, define $B_n\subseteq X$ as follows:
\[B_n=f^{-1}\left(S_n\setminus\bigcup_{k<n}S_k\right)=f^{-1}(S_n)\setminus\bigcup_{k<n}f^{-1}(S_k).\]
Since $f$ is Baire measurable, each $f^{-1}(S_n)$ has the Baire property, and since the sets with Baire property form a $\sigma$-algebra, each $B_n$ has the Baire property. Moreover if $n\neq k$ we have that $B_n\cap B_k=\varnothing$ and $X=\bigcup_{n\in\NN}B_n$.

For each $n$ let $U_n$ and $M_n$ an open set and a meager set respectively such that $B_n=U_n\mathop{\triangle}M_n$, so $U_n\cap(X\setminus M_n)=B_n\cap (X\setminus M_n)$. Let $M=\bigcup_{n\in\NN}M_n$, then $M$ is meager and for all $n\in\NN$, $B_n\cap(X\setminus M)=U_n\cap(X\setminus M)$. This means that $B_n$ is relatively open in $X\setminus M$ for all $n\in\NN$.

For every $n$ such that $B_n\setminus M$ is non-empty, choose $\xi_n\in B_n\setminus M$ and set $f_\varepsilon(x)=f(\xi_n)$ if $x\in B_n\setminus M$ (this is well-defined since the $B_n$'s are pairwise disjoint so $x$ can be at most in a single $B_n\setminus M$).

First note that $f_\varepsilon$ is continuous on $X\setminus M$, since given any $x_k\to x$ in $X\setminus M$, there is some $n$ such that $x\in B_n$ and therefore there is $j$ for which $x_k\in B_n$ for all $k\geq j$. So $f_\varepsilon(x_k)=f_\varepsilon(x)=f(\xi_n)$. We also have that $d_Y(f_\varepsilon(x),f(x))=d_Y(f(\xi_n),f(x))<\varepsilon$, as both $f(\xi_n)$ and $f(x)$ belong to the same $S_n$ which was an open ball of radius $\varepsilon$.

Repeat the construction for each $\varepsilon=\frac1k$ for $k\in\NN$ to obtain a sequence of functions $f_k$, each defined on $X\setminus N_k$, where $N_k$ is $M$ in the above construction. Then $N=\bigcup_{k\in\NN}N_k$ is a meager set and $f_k\to f$ uniformly on $X\setminus N$. Therefore $f$ is continuous as the limit of a uniformly convergent sequence of continuous functions, as wanted.
\end{proof}
Let us analyze the use of the axiom of choice in the above proof. We used $\DC$ to argue that countable unions of meager sets are meager, as well to choose each $\xi_n$ and to argue that sequential continuity is equivalent to continuity. Note that while we can appeal to Baire's theorem to conclude that $X$ is a Baire space, if $X$ is meager in itself, then the theorem is trivially true since $N=X$ satisfies the requirement.

All in all, however, it is not hard to see that indeed the proof is the same in $\ZF+\DC$. We use the cited corollary about the consistency result of Solovay and Shelah to obtain this easy corollary.

\begin{corollary}[$\DC+\BP$]
Every function $f\colon X\to Y$ between a complete and separable metric space $X$, and a separable metric space $Y$, is continuous on $X\setminus N$ for some meager set $N$.\qed
\end{corollary}

\begin{theorem}[$\DC$, Banach {\cite[Ch.~1, Th.~4]{Banach:1932}}\footnote{The proof is using Borel measurability rather than Baire measurability, but the remark after the proof explains that the argument translates mutatis mutandi to Baire measurability, as shown in this work.}]
Let $X$ be a separable Fr\'echet space and $Y$ a separable normed space. If $T\colon X\to Y$ is a linear operator which is Baire measurable, then it is continuous.
\end{theorem}
\begin{proof}
Suppose that $T\colon X\to Y$ is a linear operator and $x_n\to 0$ in $X$. By the previous lemma, there is a meager set $N\subseteq X$ such that $T\restriction X\setminus N$ is continuous. The set $M=\bigcup_{n\in\NN} -x_nN$ is meager, so there is some $z\notin M$, it follows that $x_n+z\notin M$ as well for all $n\in\NN$, since $x_n+z\to z$ and $T$ is continuous at $z$ we have that:
\begin{align*}
\lim_{n\to\infty}Tx_n
&=\left(\lim_{n\to\infty}Tx_n\right)+Tz-Tz\\
&=\left(\lim_{n\to\infty}Tx_n+Tz\right)-Tz\\
&=\left(\lim_{n\to\infty}T(x_n+z)\right)-Tz\\
&=T\left(\lim_{n\to\infty}x_n+z\right)-Tz\\
&=Tz-Tz\\
&=0.\qedhere
\end{align*}
\end{proof}
\begin{remark}
The above theorem is true in a much greater capacity, we can relax the assumptions on $Y$ (and not require separability, for example). The proof, however, is much more complicated, and we do not the broader result in this work. The proof can be found in \cite[Th.~27.45]{Schechter:1997}. However, for sake of readability, we will omit the assumptions on $Y$ for the remainder of the text.
\end{remark}
\begin{corollary}[$\DC+\BP$]
Let $X$ be any Fr\'echet space and $Y$ a normed space. If $T\colon X\to Y$ is a linear operator which is Baire measurable, then it is continuous.
\end{corollary}
\begin{proof}
Under our assumptions of $\DC+\BP$ every operator from a separable Fr\'echet space to a normed space is continuous. Suppose that $x_n\to 0$ in $X$, take $M$ to be the closure of the subspace $\operatorname{span}\{x_n\mid n\in\NN\}$, then $M$ is a separable Fr\'echet space, so $T\restriction M$ is continuous and therefore $Tx_n\to 0$ as wanted.
\end{proof}

\begin{corollary}[$\DC+\BP$]
Every solution of Cauchy's functional equation is continuous.\qed
\end{corollary}

\subsection{No finitely additive probability measures}

Another interesting consequence of $\DC+\BP$ is that there is no finitely additive probability measure on $\mathcal P(\NN)$ which vanishes on finite sets. Why is this interesting? The dual space of $L^\infty(\Omega,\Sigma,\mu)$ is composed of complex-valued finitely additive measures which are absolutely continuous relative to $\mu$. If there are no nontrivial finitely additive probability measures on $\NN$, then the dual space of $\ell^\infty$ shrinks down, and in fact becomes $\ell^1$ as we will see later.

For technical reasons it is easier to work with $2^\NN$ which is the space of all infinite binary sequences\footnote{Or the Cantor space.}, a compact metric space (and hence separable and complete). It is clear that if $\mu$ is a measure on $\mathcal P(\NN)$, we can define a function $\nu\colon 2^\NN\to[0,1]$ such that $\mu(A)=\nu(\chi_A)$. If $a\in 2^\NN$, we will write $a_n$ as the $n$th digit in the sequence $a$, and we will use $A$ to denote the set for which $a=\chi_A$, namely $A=\{n\mid a_n=1\}=a^{-1}(1)$. Finally, we will write $a'$ for $\chi_{\NN\setminus A}$, namely $a'_n=1-a_n$.

Before we proceed to the main theorem of this part, we state the following theorem by Oxtoby.\footnote{Proof of this theorem can be found in Schechter's book, \cite[Th.~20.33, p.~544]{Schechter:1997}.}

\begin{theorem}[Oxtoby's Zero-One Law]
Suppose that $T\subseteq 2^\NN$ has the Baire property and is closed under finite modifications, then $T$ is meager or $2^\NN\setminus T$ is meager.\qed
\end{theorem}

\begin{theorem}
Suppose that $\mu$ is a finitely additive probability measure on $\power(\NN)$ which vanishes on finite sets. Then $T=\{a\in 2^\NN\mid \mu(A)=0\}$ does not have the Baire property in the topological space $2^\NN$.
\end{theorem}
\begin{proof} For every $A$ let $\nu(\chi_A)=\mu(A)$, then $\nu\colon 2^\NN\to[0,1]$ and $T=\nu^{-1}(0)$.

Assume towards a contradiction that $T$ has the Baire property, by the fact that $\mu$ gives $0$ to finite sets we immediately obtain that $T$ is closed under finite modifications and thus either meager or $U=2^\NN\setminus T$ is meager. It is impossible that $U$ is meager. Suppose that it was, note that $U=\{a\mid \nu(a)>0\}$, therefore \[U'=\{a'\mid a\in U\}=\{a\mid \nu(a)<1\}.\] As $\nu(a)+\nu(a')=1$ for all $a\in 2^\NN$ it has to be the case that either $a\in U$ or $a\in U'$ for all $a\in 2^\NN$. But since $a\mapsto a'$ is a homeomorphism of $2^\NN$ onto itself, it maps meager sets to meager sets, in particular $U'$ is meager, but $2^\NN=U\cup U'$ cannot be the union of two meager sets.

It has to be the case, if so, that $T$ is meager. For every $p\in\NN$ let $Q_p$ be a closed nowhere dense set, such that $T=\bigcup_{p\in\NN}Q_p$, and $U_p=2^\NN\setminus Q_p$ is its dense open complement. We are going to derive a contradiction by showing there is some $k>0$ and pairwise disjoint sets $A_i$ for $i=1,\ldots,k$ such that $\mu(A_i)>\frac1k$ which will contradict the finite additivity of the measure.

If $a$ is a finite binary sequence of length $n$ and $V\subseteq 2^\NN$, we say that $a\in V$ if every infinite extension of $a$ belongs to $V$. If we think about this topologically, it means that the basic open set defined by $a$ is a subset of $V$. We will also write $a^\frown b$ to denote the concatenation of the finite sequences $a$ and $b$.

We will now define by induction on $p$ a strictly increasing sequence of integers $\lambda(p)$ and functions $F_p$ such that:
\begin{enumerate}
\item $F_p$ maps sequences of length $p$ to sequences of length $\lambda(p)$.
\item If $q>p$, $a$ of length $p$, and $b$ extends $a$ to length $q$, then $F_q(b)$ extends $F_p(a)$.
\item For all $a$ of length $p>0$, $F_p(a)\in U_p$.
\item For every $j=1,\ldots,\lambda(p)$ there is at most one $a$ of length $p$, such that $(F_p(a))_j=1$.
\end{enumerate}
For $p=0$ define $\lambda(p)=0$, so the empty sequence is mapped to the empty sequence. Suppose that we defined $\lambda(p)$ and $F_p$. For every $a$ of length $p$, and $i\in\{0,1\}$ let $G_p(a^\frown\tup{i})$ be $F_p(a)^\frown a^\frown\tup{i}$, which is a sequence of length $q=\lambda(p)+p+1$.

Next, enumerate all the binary sequences of length $p+1$ as $c_i$ for $i=1,\ldots,2^{p+1}$ and by induction let $\widetilde{c_i}$ be an extension of $G_p(c_i)$ such that the length of $\widetilde{c_i}$ is strictly less than the length of $\widetilde{c_{i+1}}$, and $\widetilde{c_i}\in U_{p+1}$. We can always find such sequences since $U_p$ is dense open.

Finally, set $\lambda(p+1)$ to be the length of $\widetilde{c_{2^{p+1}}}$ (the longest extension), and extend each $\widetilde{c_i}$ to this length by appending $0$'s as necessary (we will now regard $\widetilde{c_i}$ as a sequence of this length). Finally, if $b$ has length $p+1$, it is equal to some $c_i$ and define $F_{p+1}(b)=\widetilde{c_i}$. It is not hard to see that all the requirements are satisfied.

Now define a function $F\colon 2^\NN\to 2^\NN$ by letting $F(a)$ be the unique element obtained as the limit of the sequence $F_p(\tup{a_1,\ldots,a_p})$ (recall that for $p<q$ we have coherence between $F_p$ and $F_q$). By the fourth requirement $F$ has to be injective so $B=\rng(F)$ is uncountable. Moreover $B\subseteq\bigcap_{p\in\NN} U_p$ so it is disjoint from each $Q_p$ and therefore $B\cap T=\varnothing$, and for every $a,b\in B$ it follows that $A\cap B$ is finite (again, this follows from the last requirement on $F_p$).

By the uncountability of $B$ there is some $k\in\NN$ such that the set $\{b\in B\mid\nu(b)>\frac1k\}$ is uncountable, and therefore we can choose $b_1,\ldots,b_k\in B$ such that $\nu(b_i)>\frac1k$ for all $i=1,\ldots,k$. Then $\mu(B_i)>\frac1k$ for all $i=1,\ldots,k$. Now since $\mu$ is invariant under finite changes and each $B_i\cap B_j$ is finite, we can remove some finite part from all the $B_i$'s to obtain pairwise disjoint $A_i$ for $i=1,\ldots,k$ with $\mu(A_i)>\frac1k$ which is the wanted contradiction.
\end{proof}

\begin{corollary}[$\DC+\BP$]
There are no finitely additive probability measures on $\power(\NN)$ which vanish on finite sets.\qed
\end{corollary}
\section{Strange consequences}

We will investigate some of the stranger consequences of $\ZF+\DC+\BP$ on functional analysis and vector spaces in general.
\vspace{1ex}

We begin with an easy consequence of Riesz representation theorem.
\begin{theorem}
The algebraic dual of $\ell^2$ is naturally isomorphic to $\ell^2$.\qed
\end{theorem}

This theorem is interesting. One of the classical theorems of linear algebra is that if $V'$ is the \textit{algebraic} dual of $V$, then $\dim V=\dim V'$ if and only if $\dim V<\infty$. Moreover the isomorphism between $V$ and $V'$ depends on the choice of basis, whereas the isomorphism between $V$ and $V''$ is natural in the sense that it does not require a choice of basis. In our case, not only $\ell^2$ is not a finite dimensional vector space it is also naturally isomorphic to its dual (not just the double dual). Of course, $\ell^2$ has no Hamel basis in this case.

We can have an even more extreme situation, where $\dim V'$ is finite, and even that $\dim V'=0$. Before proving this, we will need to talk about Bartle integrals.
\subsection{Bartle Integrals on Sequences}

Recall that $x\in\ell^\infty$ is a function from $\NN$ to $\CC$, so $x^{-1}(\alpha)=\{n\in\NN\mid x_n=\alpha\}$ for $\alpha\in\CC$. We will define a new notion of integration for finitely additive complex-valued measures on $\power(\NN)$.

\begin{definition}
Suppose that $\mu$ is a finitely additive complex-valued measure on $\power(\NN)$ and $x\in\ell^\infty$, we define the \textit{Bartle integral}\footnote{This definition can be given in a much broader context, but we are only interested in the case of $\ell^\infty$ here. So it is better to simplify the definition to match our needs.} of $\mu$ as follows:
\[\int_\NN x\d\mu=\sum_{\alpha\in\CC}\alpha\mu(x^{-1}(\alpha)).\]
\end{definition}
\begin{proposition}
If $\mu$ is a finitely additive bounded measure on $\power(\NN)$ and $x\in\ell^\infty$, then the following inequality holds:
\[\left|\int_\NN x\d\mu\right|\leq\|x\|_\infty|\mu(\NN)|.\qed\]
\end{proposition}
If $D$ is the space of all finitely additive complex measures on $\power(\NN)$, then \[\tup{\mu,x}=\int_\NN x\d\mu\] is a continuous bilinear form on $D\times\ell^\infty$ (here the integration is Bartle). So for every $\mu\in D$ we get a $\varphi_\mu(x)=\tup{\mu,x}$, which is a continuous linear functional on $\ell^\infty$.
\begin{theorem}
The following are equivalent:\begin{enumerate}
\item There exists a finitely additive probability measure on $\power(\NN)$ which vanishes on finite sets.
\item $\ell^1\subsetneqq(\ell^\infty)^*$.
\end{enumerate}
\end{theorem}
\begin{proof}
The one direction is easy. If $\mu$ is a finitely additive probability measure on $\power(\NN)$ which vanishes on finite sets, then $\varphi_\mu$ is a continuous linear functional on $\ell^\infty$, and there is no $a\in\ell^1$ such that $\varphi_\mu(x)=\sum_n a_nx_n$, since for every $n$, the sequence $e^n=\chi_{\{n\}}$ is such that $\int_\NN e^n\d\mu=0$, and so $a_n$ would have to be $0$ for all $n$.

In the other direction, let $\varphi$ be a continuous linear functional on $\ell^\infty$ which is not defined from a sequence in $\ell^1$. The function $\mu(S)=\varphi(\chi_S)$ is a finitely additive (complex-valued) measure on $\power(\NN)$, which is not $\sigma$-additive as in that case it would be true that $\mu(\NN)=\sum_{n\in\NN}\varphi(e^n)$, which would mean that $\varphi$ is defined from a sequence in $\ell^1$. By decomposing $\mu$ to real and imaginary, and each of those to positive and negative finitely additive measures we are guaranteed that at least one of them will not be $\sigma$-additive. So without loss of generality let $\mu$ be that same measure already.

Now we claim that there is some $J\subseteq\NN$ such that $\mu(J)>\sum_{n\in J}\mu(\{n\})$, otherwise $\mu$ was $\sigma$-additive. Define the measure $\nu$ as follows:
\[\nu(S)=\frac{\mu(S\cap J)-\sum_{n\in S\cap J}\mu(\{n\})}{\mu(J)}.\]
Then $\nu$ is a finitely additive probability measure on $\power(\NN)$ which trivially vanishes on finite sets, as wanted.
\end{proof}

We infer three corollaries of this theorem, under our assumption of $\ZF+\DC+\BP$.

\begin{theorem}[V\"ath \cite{Vath:1998}]
$\ell^1$ is reflexive.
\end{theorem}
\begin{proof}
If $\BP$ holds, then there are no nontrivial finitely additive probability measures on $\power(\NN)$ and therefore every linear functional on $\ell^\infty$ comes from $\ell^1$.\footnote{Note we do not care about continuity, since we already know that under our assumptions every linear functional is continuous.}
\end{proof}
\begin{theorem}
$(\ell^\infty/c_0)^*$ is trivial.
\end{theorem}
\begin{proof}
Suppose that $\varphi$ is a linear functional on $\ell^\infty/c_0$, then there is $\widetilde\varphi\in(\ell^\infty)^*$ such that $\widetilde\varphi(x)=\varphi(x+c_0)$. However every linear functional on $\ell^\infty$ comes from $\ell^1$, therefore there is some $a\in\ell^1$ such that $\widetilde\varphi(x)=\sum_n a_nx_n$ and for every $x\in c_0$, $\widetilde\varphi(x)=0$. It follows immediately that $a_n=0$ for all $n$ and therefore $\varphi\equiv 0$ as wanted.
\end{proof}
\begin{corollary}
There is an infinite dimensional vector space whose algebraic dual is trivial.\qed
\end{corollary}

\subsubsection{Other practical consequences}

The weak topology on a Banach space and the the weak-$*$ topology on its dual are important in functional analysis. In the case that every functional is continuous we get the following interesting consequence. We continue to work in $\ZF+\DC+\BP$.

\begin{theorem}
Let $X$ be a vector space over $\CC$. Suppose that $\tau_1$ and $\tau_2$ are two topologies on $X$ which induce a Fr\'echet space structure on $X$, then $\tau_1=\tau_2$. In other words, a space $X$ can have at most one Fr\'echet topology.
\end{theorem}

\begin{proof}
Let $Tx=x$ be the identity operator, then $T$ is a linear operator from $(X,\tau_1)$ to $(X,\tau_2)$ and therefore continuous, and for the same reason $T^{-1}$ is continuous also.
\end{proof}

The example of $\ell^\infty/c_0$ is an example of a Banach space without a dual. It follows that the weak topology on $\ell^\infty/c_0$ is the indiscrete topology. But we get something else from automatic continuity. Many times we had to verify that various operators are continuous, and while the proof amounts to a possibly arduous verification of $\varepsilon$-$\delta$ arguments, under $\DC+\BP$ we only need to prove linearity to infer continuity. Notable examples include the Fourier transform, and other operators defined on the Schwartz space.\footnote{Actually, the same can be said about the space of distributions, which might not be a Fr\'echet space, but continuity is essentially equivalent to being continuous on Fr\'echet subspaces.}

\section{An open problem: Banach limits}

We will finish with a quick review of an open problem related to the Baire property and automatic continuity. For simplicity, we will work over $\RR$, rather than over $\CC$.
\begin{definition}
We say that $L\colon\ell^\infty\to\RR$ is a \textit{Banach limit} if it is a continuous linear functional satisfying:
\begin{enumerate}
\item For every $x\in\ell^\infty$, $\liminf_{n\to\infty} x_n\leq L(x)\leq\limsup_{n\to\infty} x_n$, and in particular for $x\in c$, $L(x)=\lim_{n\to\infty}x_n$.
\item $\|L\|=1$.
\item If $x\in\ell^\infty$ and $x^+$ defined as $x^+_n=x_{n+1}$, then $L(x)=L(x^+)$ (shift invariance).
\end{enumerate}
\end{definition}
The idea is that a Banach limit extends the notion of limits, but satisfies linearity whereas $\limsup$ and $\liminf$ do not.

\begin{theorem}[$\BPI$]
There exists a Banach limit.
\end{theorem}
\begin{proof}
We begin with the following lemma (whose proof does not use the axiom of choice).
\begin{lemma}
Suppose that $\mu$ is a two-valued measure on $\power(\NN)$, $\tup{t_n\mid n\in\NN}$ a bounded sequence in $\RR$, then there exists a unique $\tau\in\RR$ such that for every open $U$ such that $\xi\in U$, $\mu(\{n\mid t_n\in U\})>0$. Moreover, the function $B(\tup{t_n\mid n\in\NN})=\tau$ is a continuous linear functional and $\|B\|=1$.
\end{lemma}
\begin{proof}
Let $D=\overline{\{t_n\mid n\in\NN\}}$, then $D$ is a compact subset of $\RR$. If such $\tau$ does not exist, then for every $\tau\in D$ there is an open neighborhood $U$ such that $\mu(\{n\mid t_n\in U\})=0$, now \[\{U\mid\exists\tau\in D,\ \tau\in U,\mu(\{n\mid t_n\in U\})=0\}\] is an open cover of a compact set, so it has a finite subcover witnessed by some $\tau_1,\ldots,\tau_k$. But now $\NN=\bigcup_{i=1}^k\{n\mid t_n\in U_{\tau_i}\}$, but on the other hand this is a finite union of null sets, which is a contradiction since $\mu(\NN)>0$.

The uniqueness of $\tau$ follows from the fact $\RR$ is Hausdorff, it $\tau\neq\tau'$ then there are two disjoint open sets $U,U'$ such that $\tau\in U$ and $\tau'\in U'$. It is impossible that both $\{n\mid t_n\in U\}$ and $\{n\mid t_n\in U'\}$ have positive measure since the sets are disjoint. So it has to be the case that $\tau$ is unique. Therefore the $B$ is well-defined.

Linearity follows from the fact that $\RR$ is locally convex, namely if $B(t)=\tau$ and $\alpha\in\RR$, either $\alpha=0$ and the result is trivial, or $\alpha\neq 0$ and then whenever $U$ is an open neighborhood of $\tau$, $\alpha U$ is an open neighborhood of $\alpha\tau$ and $\{n\mid t_n\in U\}=\{n\mid\alpha t_n\in\alpha U\}$. Additivity follows similarly, using the fact that intersecting two sets of positive measure is again a set of positive measure (the measure is two-valued). Finally, $\|B\|=1$ as $|\tau|\leq\sup|t_n|=\|t\|_\infty$.
\end{proof}

One might be tempted to use $B$ as the Banach limit, but it is not shift invariant. For example, $x_n=(-1)^n$ will be mapped by $B$ to either $1$ or $-1$, but shift invariance means that $B(x)=B(x^+)$, so $B(x-x^+)=0$, but $(x-x^+)_n=(-1)^n\cdot 2$, whose value under $B$ must be $\pm2$. We fix this by considering the Ces\`aro sum operator:
\[ C(x)=y\iff y_n=\frac1n\sum_{i=1}^n x_i.\]
We claim that $C$ is a continuous linear operator, and $\|C\|=1$, then $L(x)=B(C(x))$ will be a Banach limit.

First note that $C$ is well-defined since if $\sup|x_n|<M$, then for all $n$ we have
\[ M = \frac1nnM = \frac1n\sum_{i=1}^nM\geq\frac1n\sum_{i=1}^n|x_n|\geq\left|\frac1n\sum_{i=1}^nx_n\right|=|y_n|.\]
Therefore not only $y$ is a  bounded sequence, $\|y\|_\infty\leq\|x\|_\infty$ and we get continuity (assuming linearity). It is clear that $x$ such that $x_n=1$ also satisfies $C(x)=x$ so if $C$ is linear we have $\|C\|=1$. But since linearity is trivial, the wanted conclusions follow.

We define $L(x)=B(C(x))$ where $B$ is an operator as defined in the lemma, with a two-valued measure $\mu$ which vanishes on finite sets. Such measure $\mu$ exists by $\BPI$ since the set of complement of finite sets satisfies the finite intersection property. Linearity, continuity and $\|L\|=1$ all follow from the fact that $B$ and $C$ are linear, continuous and have norm $1$ and that $L(x)=1$ where $x$ is the constant sequence $1$.

If $m=\liminf x_n$ and $M=\limsup x_n$, then $m\leq C(x)_n\leq M$ so the first requirement holds. To see that indeed $L(x)=L(x^+)$ note that $\lim C(x-x^+)_n=0$:
\[C(x)_n-C(x^+)_n=\frac1n\sum_{i=1}^nx_i-\frac1n\sum_{i=1}^nx_{i+1}=\frac1n\left(\sum_{i=1}^nx_i-x_{i+1}\right)=\frac1n(x_1-x_{n+1})\to0.\]
Therefore $L(x-x^+)=0$ so it has to be the case that $L(x)=L(x^+)$ as wanted.
\end{proof}

Having established that it is a consequence of $\AC$ that Banach limits exist (recall that $\BPI$ is just a theorem of $\AC$), the next theorem shows that we cannot quite remove all traces of choice from our construction.

\begin{theorem}[$\DC$]
If there exists a Banach limit, then $\BP$ fails.
\end{theorem}
\begin{proof}
It suffices to note that a Banach limit cannot be defined from a sequence in $\ell^1$, as for every $n$, $e_n=\chi_{\{n\}}$ must be mapped to $0$.
\end{proof}

Looking at the above theorem, one has to wonder. We saw that automatic continuity is a consequence of the Baire property holding for every set of reals. Does the existence of a Banach limit refute the statement ``Every linear functional on a Banach space is continuous''? We required that a Banach limit is continuous, so it certainly cannot be a counterexample on its own.

\begin{theorem}[$\DC$]
Every proof in the subsection ``\nameref{Subsection:Auto}'' can be modified slightly to use the assumption ``Every set of reals is Lebesgue measurable'' instead of ``Every set of reals has the Baire property''.\qed
\end{theorem}

Historically the study of automatic continuity began with Banach when he proved that a Lebesgue measurable solution to Cauchy's functional equation is continuous \cite{Banach:1920}. It was only in 1932 that he remarked that the proof is similar with the Baire property instead. And so we can rewrite the question about automatic continuity. If the existence of Banach limits implies that automatic continuity fails, it will imply that there is a non-measurable set of reals.

\begin{question}
Suppose that there exists a Banach limit, can we prove the existence of a non-measurable set of reals?
\end{question}

This can be weakened even further. Since a Banach limit is a functional not coming from $\ell^1$, it implies the existence of a finitely additive probability measure on $\power(\NN)$ which vanishes on finite sets. This leaves us with an open question of Pincus \cite{Pincus:1974}.

\begin{question}
Suppose that $\ell^1$ is not reflexive. Is there a non-measurable set of reals?
\end{question}

\chapter{The Dream Shattered: Hahn--Banach and Friends}

In this chapter we will consider the Hahn--Banach theorem and related topics, as extending functionals is much easier in the case of Banach spaces over $\RR$, and the complex case follows immediately, we will always work over $\RR$ in this chapter.
\section{Extending linear functionals}
\begin{definition}
Let $X$ be a vector space over $\RR$. We say that $p\colon X\to\RR$ is a \textit{sublinear} function if the following holds:
\begin{enumerate}
\item $p(x+y)\leq p(x)+p(y)$.
\item $p(\alpha x)=\alpha p(x)$ for all $\alpha\geq 0$.
\end{enumerate}
\end{definition}

\begin{theorem}[Finite Extension Theorem]
Suppose that $X$ is a vector space over $\RR$ and $M\subseteq X$ is a subspace of $X$. Assume that $f_0\colon M\to\RR$ is a linear functional and $p\colon X\to\RR$ is a sublinear function such that $f_0(x)\leq p(x)$ for all $x\in M$. If $y\in X$ is any vector, then there is an extension $f$ of $f_0$ to $N=\spnn(M\cup\{y\})$ such that $f(x)\leq p(x)$ for all $x\in N$.
\end{theorem}
\begin{proof}
If we can extend $f$ to $N$, then $f(m+\alpha y)=f(m)+\alpha f(y)$ for all $m+\alpha y\in N$, but we also want to satisfy the property that for $\alpha>0$, \[f(m+\alpha y)=f(m)+\alpha f(y)\leq p(m+\alpha y).\]

Suppose that we managed to find such extension of $f$. Let $\alpha>0$, then $f(m+\alpha y)\leq p(m+\alpha y)$, but from the sublinearity of $p$ and linearity of $f$ we get:
\[\alpha\left(f\left(\frac m\alpha\right)+f(y)\right)\leq\alpha p\left(\frac m\alpha+y)\right).\]
We can divide both sides by $\alpha$ and obtain that $f(y)\leq p\left(\frac m\alpha+y\right)-f\left(\frac m\alpha\right)$. If $\alpha<0$ we can use $-\alpha$ and get that \[f(y)\geq f\left(\frac m\alpha\right)-p\left(-\left(\frac m\alpha+y\right)\right),\]
so it is always the case that $p(m_1+y)-f(m_1)\geq f(m_2)-p(m_2+y)$.

Finally, for $m_1,m_2\in M$ we already know that
\begin{align*}
 f(m_1)+f(m_2)&=f(m_1+m_2)\\ &\leq p(m_1+m_2)\\&=p(m_1+y+m_2-y)\\&\leq p(m_1+y)-p(m_2-y).
\end{align*}
And so we get the condition that
\[\inf_{m_1\in M}\biggl(p(m_1+y)-f(m_1))\biggr)\geq f(y)\geq\sup_{m_2\in M}\biggl(f(m_2)-p(m_2-y)\biggr).\]

So we can choose any $c$ in the closed interval and extend $f$ accordingly to $N$.
\end{proof}
\begin{theorem}[$\BPI$, Hahn--Banach]
Under the assumptions of the Finite Extension Theorem there exists an extension of $f$ to a linear functional on $X$ satisfying $f(x)\leq p(x)$ for all $x\in X$.
\end{theorem}
The initial instinct is to appeal to Zorn's lemma. However we want to prove it from a strictly weaker choice principle, $\BPI$. So a slightly different approach is needed here.
\begin{proof}
Let $D$ be the set of all possible extensions of $f_0$ to subspaces of $X$, then by the Finite Extension Theorem, for every $x\in X$ the set $D_x=\{g\in D\mid x\in\dom g\}$ is non-empty. Moreover by a simple induction, the family $\{D_x\mid x\in X\}$ has the finite intersection property: if $x_1,\ldots,x_n\in X$ then we can first extend $f$ to $f_1$ on $\spnn(M\cup\{x_1\})$, then from $f_i$ to $f_{i+1}$ on $\spnn(M\cup\{x_1,\ldots,x_{i+1}\})$ by applying the Finite Extension Theorem repeatedly.

From $\BPI$ we get that there is a binary measure $\mu$ on $\power(D)$ such that $\mu(D_x)=1$ for all $x\in X$. For every $x\in X$ the set $\{g(x)\mid g\in D_x\}$ is bounded between $-|p(x)|$ and $|p(x)|$. Using an argument similar to the existence of Banach limits on $\ell^\infty$ we can prove that there is a unique $c\in\RR$ such that for every $U$ open neighborhood of $c$, $\{g\in D_x\mid g(x)\in U\}$ has measure $1$. Therefore we define $f(x)$ to be that $c$.

Then $f$ extends the functional $f_0$, and it is dominated by $p$. To see that $f(x+\alpha y)=f(x)+\alpha f(y)$ note that the set of extensions which are defined on both $x$ and $y$ has measure $1$, and the proof now goes as in the Banach limit case again.
\end{proof}

It is not a coincidence that the proof is very similar to the proof of existence of Banach limits. Schechter formulates the Hahn--Banach theorem in several equivalent versions, in the first one is in fact the existence of generalized Banach limits (replacing $\NN$ by directed sets and nets).

As the Hahn--Banach theorem cannot be proved in $\ZF$ (as the next theorem shows), and it follows from $\ZF+\BPI$, one might ask if the Hahn--Banach theorem is in fact equivalent to $\BPI$. This was the prevailing conjecture in the 1960s (see Luxemberg \cite{Lux:1962}), but Pincus disproved it in 1974 \cite{Pincus:1974}.

We next prove the following theorem, which should be fairly trivial at this point.
\begin{theorem}
Hahn--Banach cannot be proved in $\ZF$.
\end{theorem}
\begin{proof}
Pick any non-zero $x\in\ell^\infty/c_0$ and consider the functional $f_0(\alpha x)=\alpha$ as defined on $\spnn\{x\}$. By the Hahn--Banach theorem this can be extended to a linear functional on $\ell^\infty/c_0$ which is non-zero. Therefore there exists a set of real numbers without the Baire property.
\end{proof}

The above proof, albeit simple, shows that sometimes one linear functional is enough to deduce to much more. This brings us to the following theorem by Luxemburg and V\"ath \cite{Lux-Vath:2001} which shows that just like bedbugs, fleas and other pests, if there is one continuous linear functional on each Banach space, then there are many.

\begin{theorem}[$\DC$]
Hahn--Banach is equivalent to the statement that every nontrivial Banach space has a continuous non-zero functional.
\end{theorem}
\begin{proof}
One direction is trivial. If $X$ is a Banach space, pick any $x_0\neq 0$ such that $\|x_0\|_X=1$ and look at $f(\alpha x_0)=\alpha$ as a functional on the one-dimensional subspace $\spnn\{x_0\}$. By the Hahn--Banach theorem this extends to a functional bounded by the norm, i.e.\ continuous.

Now suppose that $X$ is a vector space over $\RR$, $M\subseteq X$ a subspace with $f_0\colon M\to\RR$ a linear functional bounded by a sublinear function $p\colon X\to\RR$. Let $D$ be the set of all extensions of $f_0$ to subspaces of $X$ which are still dominated by $p$. For every finite $S\subseteq X$ let $D_S$ be the set of extensions which are defined on $M\cup S$.

Consider the space $Y=\ell^\infty(D)$ whose elements are bounded functions $y\colon D\to\RR$, endowed with the $\sup$ norm this is a Banach space. Next, consider the subspace $Y_0$ defined as:
\[ Y_0=\{y\in Y\mid\exists S\subseteq X\text{ finite }, y\restriction D_S\equiv 0\}.\]
By the Finite Extension Theorem, for every finite $S$ the set $D_S$ is non-empty, so there is $y\in Y_0$ such that $y\restriction D_S\equiv 0$. Let $Z=Y/\overline{Y_0}$ be the quotient space, which is a Banach space as $\overline{Y_0}$ is closed. For $y\in Y$ we will write $[y]$ to be the equivalence class of $y$ in $Z$.

First note that $e_0\equiv 1$ is such that $d(Y_0,e_0)=1$ and therefore $Z$ is not a trivial Banach space. It also follows that $e=[e_0]$ is such that $\|e\|_Z=1$.

Let $G\colon X\to Z$ be defined as follows, for every $x\in X$ look at the evaluation functional, $y(f)=f(x)$ and extended it to be bounded somehow on the set of $f$'s for which $x\notin\dom f$. Note that if $y_1,y_2$ are two extensions, then for every $f\in D_{\{x\}}$ it holds that $y_1(f)=y_2(f)=f(x)$ and therefore $(y_1-y_2)(f)=0$. So $[y_1]=[y_2]$, and therefore defining $G(x)=[y]$ where $y$ is any bounded extension like that, is a well-defined function, and clearly it is linear.

By our assumption $Z$ has a nontrivial continuous functional $L$, define $F_0(x)=L(G(x))$, then $F_0$ is a linear functional on $X$. We note that the following holds for every $x\in X$,
\[ F_0(x)=L(G(x))=L([f\mapsto f(x)])\leq L([p(x) e_0])=p(x) L(e).\]
In particular $L(e)\neq 0$ so by normalizing we may assume $L(e)=1$, so $F_0$ is dominated by $p$. It is clear that $F_0$ extends $f_0$, since for every $m\in M$, $f(m)=f_0(m)$ for every $f\in D$ and therefore $G(m)=[f\mapsto f_0(m)]$.
\end{proof}
\begin{remark}
The above proof seems not to be using the assumption of $\DC$. But in fact that assumption is used to ensure that $Z$ is a Banach space. If we rephrase the equivalence to ``Every normed space has a nontrivial bounded functional'' then the proof is the same in $\ZF$, and we can infer the Hahn--Banach theorem.
\end{remark}
\begin{corollary}If the Hahn--Banach theorem fails, then there exists a Banach space with a trivial (continuous) dual.\qed
\end{corollary}

We will also mention that the ``infamous'' Banach--Tarski paradox which asserts that the unit ball of $\RR^3$ can be partitioned into finitely many parts and by rotations and shifts reassembled into two copies of the unit ball, can be proved from the Hahn--Banach theorem. This was shown by Pawlikowski in 1991 \cite{Paw:1991}.

\subsection{What can we do without choice?}

In the case of Baire's theorem we saw that despite the necessity of the axiom of choice for the general case, if we are willing to restrict ourselves to separable spaces then we can still prove the theorem. This is true in the Hahn--Banach case as well.

\begin{theorem}
Suppose that $X$ is a separable Banach space, $M_0$ a subspace of $X$ and $f_0$ a linear functional on $M_0$ such that $f_0(m)\leq c\|m\|_X$ for some constant $c>0$. Then there is an extension of $f_0$ to $f\colon X\to\RR$ such that $f$ is continuous.
\end{theorem}
Before we prove this, let us remark that while many of the spaces of interest are separable, $\ell^\infty$ is not separable. And as we have seen before, it is possible that Hahn--Banach fails for $\ell^\infty$.

\begin{proof}
Let $\{x_n\mid n\in\NN\}$ be a fixed countable dense subset, and define by induction $M_{n+1}$ to be $\spnn(M_n\cup\{x_n\})$, where we extend $f_n$ to $f_{n+1}$ by defining $f_{n+1}(x_n)$ to be the maximum possible value as in the proof of the Finite Extension Theorem, as the values in the proof are taken from a canonically defined closed interval, we have a well-defined choice of value for the extension.

Let $M=\bigcup_{n\in\NN} M_n$, then $M$ is a dense subspace of $X$ and $f$ defined as $f(x)=f_n(x)$ if $x\in M_n$ (as $f_n$ extends $f_m$ for $m\leq n$, this is well-defined). Now since $f$ is dominated by $p$ we can extend it uniquely to a continuous functional on $X$.
\end{proof}

The above proof might seem suspicious. If it works, why doesn't the general proof work? And if the general proof fails, are we sure that we did not use the axiom of choice in \textit{this} proof? The answer is that the fact we can choose canonically a value for $f(x_n)$ means that we extended $f_0$ to the countable subset in a well-defined manner. In the general proof we did not have a dense countable subset either, even if each one-step extension can be done canonically, we would not have a concrete set to use induction for the extension.\footnote{Generally in choiceless set theory most proofs that work in $\ZFC$ can be repeated when the sets involved are countable or well-orderable in general. So when we assume that we have a ``generating set'' in some sense (like a dense subset, when continuous functions are involved) which is well-orderable, we can sometimes use that to our advantage and avoid the axiom of choice altogether.}

\section{The Banach--Alaoglu and Krein--Milman theorems}

Remember that if $X$ is a Banach space, its (continuous) dual space $X^*$ has three naturally arising topologies defined: The strong topology, the weak topology and the weak-$*$ topology. The following is an important theorem of functional analysis,

\begin{theorem}[$\AC$, Banach--Alaoglu]
If $X$ is a normed space, then the closed unit ball of $X^*$ is compact in the weak-$*$ topology.\qed
\end{theorem}

As we have seen, $X^*$ might be trivial in which case it is certainly true that the unit ball is compact. But this does not mean that the theorem is provable in $\ZF$ without any additional assumption. Indeed, the following theorem shows that this is not the case.

\begin{theorem}
The following are equivalent:
\begin{enumerate}
\item The Banach--Alaoglu Theorem.
\item $\BPI$.
\end{enumerate}
\end{theorem}
\begin{proof}
Assume that the Banach--Alaoglu theorem holds. Let $S$ be a non-empty set, and suppose that $\cF$ is a collection of subsets of $S$ with the finite intersection property. Denote by $X$ the space $\ell^\infty(S)$, then by the Banach--Alaoglu theorem the unit ball $B$ of $X^*$ is weak-$*$ compact.

For every $\varphi\in X^*$ we write $\mu_\varphi$ to be the function on $\mathcal P(S)$ defined by $\mu_\varphi(A)=\varphi(\chi_A)$. And since $\varphi(\chi_\varnothing)=0$, $\mu_\varphi$ is a finitely additive signed measure.\footnote{Last time we did something like this we got complex-valued measures, but we worked over $\CC$ and now we work over $\RR$ so we only get signed measures.} Then the following set is a closed subset of $B$, \[U=\{\varphi\in B\mid\mu_\varphi(S)=1,\forall A,\mu_\varphi(A)\in\{0,1\}\},\] and therefore it is compact in the weak-$*$ topology.

To see that $U$ is closed, suppose that $\varphi\notin U$ and let $A\subseteq S$ such that $1>\varphi(\chi_A)=\varepsilon>0$. Let $0<\delta<\min\{\varepsilon,1-\varepsilon\}$, then the open set $V=\varphi+U^\delta_A=\{\psi\mid |\psi(\chi_A)-\varphi(\chi_A)|<\delta\}$ is a weak-$*$ open neighborhood of $\varphi$ and if $\psi\in U\cap V$ then $\psi(A)\in\{0,1\}$ in which case $|\psi(A)-\varphi(A)|>\delta$ which is a contradiction.

Next, for every $a\in S$ the linear functional $\delta_a(f)=f(a)$ is a continuous linear functional and it is clear that $\delta_a\in U$ for all $a\in S$. For every $A\in\cF$ denote by $D_A=\overline{\{\delta_a\mid a\in A\}}$, this is a closed subset of $U$ and $\{D_A\mid A\in\cF\}$ has the finite intersection property, by the fact $\cF$ has the finite intersection property.

By compactness of $U$ there is some $\varphi\in\bigcap_{A\in\cF}D_A$, then $\mu_\varphi$ is a binary measure on $\mathcal P(S)$ as wanted.

In the other direction we note that the usual proof of the Banach--Alaoglu works out of the box, simply by noting that $\BPI$ is equivalent to the Tychonoff theorem restricted to Hausdorff spaces.\footnote{This equivalence is easy to obtain when talking about ultrafilters and characterizations of compactness by ultrafilters, but these lie beyond the intended scope of this work.}
\end{proof}

On its own the Banach--Alaoglu theorem is very useful. Compact sets are in some deep sense of the word almost-finite, and it means that many times, using them is almost like working with finite sets. One of the most useful theorems about compact sets in locally convex vector spaces is due to Krein and Milman, about the existence of extremal points. And this will be the topic of the final part of the work.

\subsection{Krein--Milman theorems}

\begin{theorem}[$\AC$, Krein--Milman]
If $X$ is a locally convex vector space, $K$ is a compact and convex set then $K$ is the closure of the convex hull of the extremal points of $K$.\qed
\end{theorem}

The proof of this theorem goes through the following lemma, whose proof relies heavily on Zorn's lemma.

\begin{theorem}[$\AC$, Krein--Milman Lemma]
In a locally convex vector space, every non-empty compact and convex set has an extremal point.\qed
\end{theorem}

Since our interest lie in the necessity of the axiom of choice, and thus Zorn's lemma, in functional analysis, it is natural to ask whether or not we can prove this lemma without the axiom of choice. All the following proofs are due to Bell and Fremlin in \cite{Bell-Fremlin:1972}, unless stated otherwise.

\begin{theorem}\label{KM+BA=AC}
Assume that Krein--Milman Lemma and the Banach--Alaoglu theorem are true, then the axiom of choice is true.
\end{theorem}

To prove this theorem, we will first prove that a different statement implies the axiom of choice, and that statement follows from the assumptions of Krein--Milman Lemma and Banach--Alaoglu trivially.

\begin{theorem}\label{Extremal Choice}
Suppose that for every normed space $X$, the unit ball of $X^*$ has an extremal point, then the axiom of choice is true.
\end{theorem}
\begin{proof}
Let $\{A_i\mid i\in I\}$ be a family of non-empty sets. Without loss of generality $A_i\cap A_j=\varnothing$ for $i\neq j$, otherwise replace $A_i$ by $\{i\}\times A_i$ for all $i\in I$. Let $A=\bigcup_{i\in I}A_i$ and consider the following space,
\[E=\left\{x\colon A\to\RR\mathrel{}\middle|\mathrel{}\forall\varepsilon>0\{t\in A\mid |x(t)|>\varepsilon\}\text{ is finite, and }\sum_{i\in I}\sup_{t\in A_i}|x(t)|<\infty\right\}.\]
It is not hard to see why $E$ is a vector space over $\RR$, and that $\|x\|_E=\sum_{i\in I}\sup_{t\in A_i}|x(t)|$ is a norm on $E$. Consider the dual space $E^*$, it can be represented as follows:
\[E^*=\left\{x\colon A\to\RR\mathrel{}\middle|\mathrel{}\sup_{i\in I}\sum_{t\in A_i}|x(t)|<\infty\right\},\]
along with the norm $\|x\|_*=\sup_{i\in I}\sum_{t\in A_i}|x(t)|$. By the assumption the unit ball of $E^*$ has an extremal point, $e$. We claim that for each $i\in I$ there exists a unique $t\in A_i$ such that $e(t)\neq 0$, and therefore $e$ defines a choice function from $\{A_i\mid i\in I\}$.

First to see there is some $t\in A_i$ for which $e(t)\neq 0$. Assume towards contradiction that $i$ was a counterexample, and fix some $v\in A_i$. Define $x,y\in X^*$ as follows:
\[ x(v)=1,\quad y(v)=-1,\quad x(t)=y(t)=e(t)\ \text{ for every }t\neq v.\]
It is not hard to see that $x,y$ lie on the unit ball of $E^*$ and that $e=\frac12(x+y)$ but $e\neq x,y$ which contradicts the assumption that $e$ was an extremal point.

Now assume towards contradiction that we have $u,v\in A_i$ such that $e(u)$ and $e(v)$ are both non-zero. We define $x,y\in E^*$ again:
\begin{align*}
x(u) &= e(u)(1+|e(v)|)\\
x(v) &= e(v)(1-|e(u)|)\\
y(u) &= e(u)(1-|e(v)|)\\
y(v) &= e(v)(1+|e(u)|)\\
x(t) &= y(t) = e(t) \quad\text{ for every }t\neq u,v.
\end{align*}
It is again easy to see that $e=\frac12(x+y)$, but $e\neq x,y$ and as before this is a contradiction to the assumption that $e$ is an extremal point. Therefore $e(t)\neq 0$ for a unique $t\in A_i$ for all $i\in I$, and so it defines a choice function as wanted.
\end{proof}

\begin{proof}[Proof of \autoref{KM+BA=AC}]
If $X$ is a normed space, then by the Banach--Alaoglu theorem the unit ball of $X^*$ is compact in the weak-$*$ topology and by Krein--Milman Lemma it has an extremal point. Therefore the assumption of \autoref{Extremal Choice} is satisfied, so the axiom of choice holds.
\end{proof}

\begin{remark}
Note that we can also work with Krein--Milman theorem and Banach--Alaoglu instead of the lemma. If the unit ball of $X^*$ is compact, it has to be the closure of the convex hull of its extremal points, and in particular it will have an extremal point.
\end{remark}

\begin{corollary}
If the axiom of choice fails but $\BPI$ holds, then the Krein--Milman theorem fails.
\end{corollary}
\begin{proof}
If the axiom of choice fails then there is a normed space $X$ such that the unit ball of $X^*$ has no extremal points. Assuming that $\BPI$ holds, the unit ball is compact (in the weak-$*$ topology) but without extremal points it cannot be the closure of the convex hull of its extremal points.
\end{proof}

We have remarked that $\BPI$ is stronger than the Hahn--Banach theorem. Perhaps we can weaken the use of $\BPI$ to that of Hahn--Banach? The answer is almost positive. Namely, we can weaken $\BPI$ to Hahn--Banach but we have to pay by strengthening Krein--Milman Lemma slightly.

\begin{definition}
Let $X$ be a normed space, we say that $A\subseteq X$ is \textit{convex-compact}\footnote{Bell and Fremlin call this property \textit{quasicompact}, but we chose to follow Luxemburg in \cite{Lux:1969} as the name fits slightly better.} if whenever $\cF$ is a family of closed convex sets satisfies that $\{A\cap F\mid F\in\cF\}$ has the finite intersection property, then $A\cap\bigcap_{F\in\cF} F\neq\varnothing$.
\end{definition}

Note that this definition extends the usual definition of compactness.

\begin{theorem}[$\AC$, Krein--Milman Strong Lemma]
Every convex and convex-compact in a locally convex space has an extremal point.\qed
\end{theorem}

We will use the following theorem, due to Luxemburg \cite{Lux:1969}, whose proof is beyond the scope of this work.
\begin{theorem}
The Hahn--Banach theorem is equivalent to the following statement: if $X$ is a normed space, the unit ball of $X^*$ is convex-compact in the weak-$*$ topology.\qed
\end{theorem}

These two theorems imply, as with \autoref{KM+BA=AC}, the following corollary
\begin{theorem}
Assume that the Krein--Milman Strong Lemma and the Hahn--Banach theorem hold, then the axiom of choice is true.\qed
\end{theorem}

We finish by pointing out that while we know that $\BPI$ is strictly stronger than Hahn--Banach, in terms of choice needed for the proof, it is not known whether Krein--Milman Strong Lemma is equivalent to the weak Lemma without the axiom of choice.
\bibliographystyle{amsplain}

\begin{thebibliography}{10}

\bibitem{Banach:1920}
Stefan Banach, \emph{Sur l'{\'e}quation fonctionnelle f (x+ y)= f (x)+ f (y)},
  Fund. Math. \textbf{1} (1920), no.~1, 123--124.

\bibitem{Banach:1932}
\bysame, \emph{Th\'eorie des op\'erations lin\'eaires}, \'Editions Jacques
  Gabay, Sceaux, 1993, Reprint of the 1932 original. \MR{1357166 (97d:01035)}

\bibitem{Bell-Fremlin:1972}
J.~L. Bell and D.~H. Fremlin, \emph{A geometric form of the axiom of choice},
  Fund. Math. \textbf{77} (1972), no.~2, 167--170. \MR{0327523 (48 \#5865)}

\bibitem{Blair:1977}
Charles~E. Blair, \emph{The {B}aire category theorem implies the principle of
  dependent choices}, Bull. Acad. Polon. Sci. S\'er. Sci. Math. Astronom. Phys.
  \textbf{25} (1977), no.~10, 933--934. \MR{0469765 (57 \#9546)}

\bibitem{Brunner:1987}
Norbert Brunner, \emph{Garnir's dream spaces with {H}amel bases}, Arch. Math.
  Logik Grundlag. \textbf{26} (1987), no.~3-4, 123--126. \MR{905950
  (88g:03066)}

\bibitem{Cohen:1963}
Paul Cohen, \emph{The independence of the continuum hypothesis}, Proc. Nat.
  Acad. Sci. U.S.A. \textbf{50} (1963), 1143--1148. \MR{0157890 (28 \#1118)}

\bibitem{Cohen:1964}
Paul~J. Cohen, \emph{The independence of the continuum hypothesis. {II}}, Proc.
  Nat. Acad. Sci. U.S.A. \textbf{51} (1964), 105--110. \MR{0159745 (28 \#2962)}

\bibitem{Enderton:Elements}
Herbert~B. Enderton, \emph{Elements of set theory}, Academic Press [Harcourt
  Brace Jovanovich, Publishers], New York-London, 1977. \MR{0439636 (55
  \#12522)}

\bibitem{Fellhauer:2015}
Adrian F.~D. Fellhauer, \emph{On the relation of three theorems of analysis to
  the axiom of choice}, J. Log. Anal. \textbf{9} (2017), Paper No. 1, 23,
  arXiv/1509.01078. \MR{3646652}

\bibitem{Garnir:1973}
H.~G. Garnir, \emph{Solovay's axiom and functional analysis}, Proceedings of
  the {S}ymposium on {F}unctional {A}nalysis ({I}stanbul, 1973), Math. Res.
  Inst., Istanbul, 1974, pp.~57--68. Publ. Math. Res. Inst.--Istanbul, No. 1.
  \MR{0477688 (57 \#17203b)}

\bibitem{Godel:AC}
Kurt G{\"o}del, \emph{The {C}onsistency of the {C}ontinuum {H}ypothesis},
  Annals of Mathematics Studies, no. 3, Princeton University Press, Princeton,
  N. J., 1940. \MR{0002514 (2,66c)}

\bibitem{Halbeisen:CombST}
Lorenz~J. Halbeisen, \emph{Combinatorial set theory}, Springer Monographs in
  Mathematics, Springer, London, 2012, With a gentle introduction to forcing.
  \MR{3025440}

\bibitem{Herrlich:AC}
Horst Herrlich, \emph{Axiom of choice}, Lecture Notes in Mathematics, vol.
  1876, Springer-Verlag, Berlin, 2006. \MR{2243715 (2007f:03069)}

\bibitem{HowardRubin:Consequences}
Paul Howard and Jean~E. Rubin, \emph{Consequences of the axiom of choice},
  Mathematical Surveys and Monographs, vol.~59, American Mathematical Society,
  Providence, RI, 1998, With 1 IBM-PC floppy disk (3.5 inch; WD). \MR{1637107
  (99h:03026)}

\bibitem{Jech:ST}
Thomas Jech, \emph{Set theory}, Springer Monographs in Mathematics,
  Springer-Verlag, Berlin, 2003, The third millennium edition, revised and
  expanded. \MR{1940513 (2004g:03071)}

\bibitem{Jech:AC}
Thomas~J. Jech, \emph{The axiom of choice}, North-Holland Publishing Co.,
  Amsterdam-London; Amercan Elsevier Publishing Co., Inc., New York, 1973,
  Studies in Logic and the Foundations of Mathematics, Vol. 75. \MR{0396271 (53
  \#139)}

\bibitem{Kunen:2011Book}
Kenneth Kunen, \emph{Set theory}, Studies in Logic (London), vol.~34, College
  Publications, London, 2011. \MR{2905394}

\bibitem{Los-RyllN}
J.~{\L}o{\'s} and C.~Ryll-Nardzewski, \emph{On the application of {T}ychonoff's
  theorem in mathematical proofs}, Fund. Math. \textbf{38} (1951), 233--237.
  \MR{0048795 (14,70h)}

\bibitem{Lux:1962}
W.~A.~J. Luxemburg, \emph{Two applications of the method of construction by
  ultrapowers to anaylsis}, Bull. Amer. Math. Soc. \textbf{68} (1962),
  416--419. \MR{0140417 (25 \#3837)}

\bibitem{Lux:1969}
\bysame, \emph{Reduced powers of the real number system and equivalents of the
  {H}ahn-{B}anach extension theorem}, Applications of {M}odel {T}heory to
  {A}lgebra, {A}nalysis, and {P}robability ({I}nternat. {S}ympos., {P}asadena,
  {C}alif., 1967), Holt, Rinehart and Winston, New York, 1969, pp.~123--137.
  \MR{0237327 (38 \#5616)}

\bibitem{Lux-Vath:2001}
W.~A.~J. Luxemburg and M.~V{\"a}th, \emph{The existence of non-trivial bounded
  functionals implies the {H}ahn-{B}anach extension theorem}, Z. Anal.
  Anwendungen \textbf{20} (2001), no.~2, 267--279. \MR{1846601 (2002g:46005)}

\bibitem{Moore:AC}
Gregory~H. Moore, \emph{Zermelo's axiom of choice}, Studies in the History of
  Mathematics and Physical Sciences, vol.~8, Springer-Verlag, New York, 1982,
  Its origins, development, and influence. \MR{679315 (85b:01036)}

\bibitem{Paw:1991}
Janusz Pawlikowski, \emph{The {H}ahn-{B}anach theorem implies the
  {B}anach-{T}arski paradox}, Fund. Math. \textbf{138} (1991), no.~1, 21--22.
  \MR{1122274 (92m:03078)}

\bibitem{Pincus:1974}
David Pincus, \emph{The strength of the {H}ahn-{B}anach theorem}, Victoria
  {S}ymposium on {N}onstandard {A}nalysis ({U}niv. {V}ictoria, {V}ictoria,
  {B}.{C}., 1972), Springer, Berlin, 1974, pp.~203--248. Lecture Notes in
  Math., Vol. 369. \MR{0476512 (57 \#16072)}

\bibitem{Pincus-Solovay:1977}
David Pincus and Robert~M. Solovay, \emph{Definability of measures and
  ultrafilters}, J. Symbolic Logic \textbf{42} (1977), no.~2, 179--190.
  \MR{0480028 (58 \#227)}

\bibitem{RubinRubin:Equivalents}
Herman Rubin and Jean~E. Rubin, \emph{Equivalents of the axiom of choice.
  {II}}, Studies in Logic and the Foundations of Mathematics, vol. 116,
  North-Holland Publishing Co., Amsterdam, 1985. \MR{798475 (87c:04004)}

\bibitem{Schechter:1997}
Eric Schechter, \emph{Handbook of analysis and its foundations}, Academic
  Press, Inc., San Diego, CA, 1997. \MR{1417259 (98b:00009)}

\bibitem{Shelah:1984}
Saharon Shelah, \emph{Can you take {S}olovay's inaccessible away?}, Israel J.
  Math. \textbf{48} (1984), no.~1, 1--47. \MR{768264 (86g:03082a)}

\bibitem{Solovay:1970}
Robert~M. Solovay, \emph{A model of set-theory in which every set of reals is
  {L}ebesgue measurable}, Ann. of Math. (2) \textbf{92} (1970), 1--56.
  \MR{0265151 (42 \#64)}

\bibitem{Vath:1998}
Martin V{\"a}th, \emph{The dual space of {$L_\infty$} is {$L_1$}}, Indag. Math.
  (N.S.) \textbf{9} (1998), no.~4, 619--625. \MR{1691998 (2000h:46037)}

\bibitem{Wright:1973}
J.~D.~Maitland Wright, \emph{All operators on a {H}ilbert space are bounded},
  Bull. Amer. Math. Soc. \textbf{79} (1973), 1247--1250. \MR{0328649 (48
  \#6991)}

\end{thebibliography}
\providecommand{\bysame}{\leavevmode\hbox to3em{\hrulefill}\thinspace}
\providecommand{\MR}{\relax\ifhmode\unskip\space\fi MR }
\providecommand{\MRhref}[2]{%
  \href{http://www.ams.org/mathscinet-getitem?mr=#1}{#2}
}
\providecommand{\href}[2]{#2}

\end{document}